\documentclass[a4paper,11pt]{article}

\usepackage{amsthm,amsmath,amsfonts,amssymb}
\usepackage[numbers]{natbib}
\usepackage{graphicx}
\usepackage{color}   
\usepackage[usenames,dvipsnames]{xcolor}
\usepackage{authblk}

\usepackage{hyperref}
\hypersetup{pdfpagemode=FullScreen,  
backref,
colorlinks=true,
citecolor=Bittersweet,
linkcolor=Bittersweet,
urlcolor=Bittersweet
}

\numberwithin{equation}{section}
\theoremstyle{plain}
\newtheorem{lemma}{Lemma}[section]
\newtheorem{theorem}{Theorem}[section]

\newcommand{\cqfd}{\hfill $\square$}

\begin{document}

\title{On the Consistency of Incomplete \mbox{U-statistics}
 under Infinite Second-order Moments}
%
%
%

\author[1]{Alexander D\"{u}rre}
\author[1,2]{Davy Paindaveine}
\affil[1]{Universit\'{e} libre de Bruxelles}
\affil[2]{Universit\'{e} Toulouse 1 Capitole}
\date{}                     

\setcounter{Maxaffil}{0}
\renewcommand\Affilfont{\itshape\small}

\maketitle


\abstract{   
We derive a consistency result, in the $L_1$-sense, for incomplete \mbox{U-statistics} in the non-standard case where the kernel at hand has infinite second-order moments. Assuming that the kernel has finite moments of order~$p(\geq 1)$, we obtain a bound on the $L_1$ distance between the incomplete \mbox{U-statistic} and its Dirac weak limit, which allows us to obtain, for any fixed~$p$, an upper bound on the consistency rate. Our results hold for most classical sampling schemes that are used to obtain incomplete \mbox{U-statistics}. 
}


\section{Introduction}
\label{secIntro}
\vspace{2mm}

The concept of \mbox{U-statistics}, that was introduced in~\cite{Hoeff1948}, has met a tremendous success in probability and statistics. In this paper, we consider throughout random variables~$X_1,X_2,\ldots$ that are independent copies of a random variable taking values in a generic space~$\mathcal{X}$ (that does not need to be~$\mathbb{R}^d$), and we fix a "kernel" function~$h:\mathcal{X}^\ell\to \mathbb{R}$ 
that is invariant under permutations of its arguments and that satisfies
\begin{equation}
\label{firstmomentassumption}
{\rm E}[
|h(X_1,\ldots,X_\ell)|
]
<
\infty
\end{equation}
(throughout the paper, $\ell$ is a fixed positive integer).
The resulting \mbox{U-statistic}, namely 
\begin{equation}
\label{Un}
U_{n}
:=
\frac{1}{{n \choose \ell}}
\sum_{1\leq i_1<\cdots<i_{\ell}\leq n}
h(X_{i_1},\ldots,X_{i_\ell})
,
\end{equation}
is then an unbiased estimator for
$\theta
:=
{\rm E}[
h(X_1,\ldots,X_\ell)
]
$.
Under these conditions, $U_n$ converges almost surely to~$\theta$ as~$n\to\infty$ (see, e.g., Theorem~A on page~190 from \cite{Ser1980}), which generalizes the classical strong law of large numbers that is obtained for~$\ell=1$. The asymptotic distribution theory for \mbox{U-statistics} was obtained under the second-order moment assumption stating that~${\rm E}[h^2(X_1,\ldots,X_\ell)]<\infty$; see, e.g., Section~5.5 from~\cite{Ser1980} or Chapter~12 from~\cite{van1998}. In the standard case of "non-degenerate" \mbox{U-statistics}, $\sqrt{n}(U_n-\theta)$ is asymptotically normal with mean zero, so that $U_n-\theta=O_P(1/\sqrt{n})$. For degenerate \mbox{U-statistics}, the consistency rate improves from~$1/\sqrt{n}$ to at least~$1/n$.

Despite the many nice properties of \mbox{U-statistics}, practical applications are often jeopardized by the fact that it may be a heavy computational burden to evaluate the $O(n^\ell)$ terms in~$U_n$. This was the motivation to introduce the \emph{incomplete \mbox{U-statistics}} 
\begin{equation} 
\label{UnN}
U_{n,N} 
:=
\frac{1}{N}
\sum_{1\leq i_1<\cdots<i_{\ell}\leq n} 
\alpha_n(i_1,\ldots,i_{\ell})
h(X_{i_1},\ldots,X_{i_\ell}) 
,
\end{equation}
where the random variables~$\alpha_n(i_1,\ldots,i_{\ell})$, with values in~$\{0,1\}$, take value one if and only if the term with index~$(i_1,\ldots,i_{\ell})$ was selected when drawing~$N$ terms without replacement out of the~${n\choose \ell}$ terms in~$U_n$. 
 Incomplete \mbox{U-statistics} were first defined in \cite{Blom1976} and were studied, among others, in~\cite{BK1978}, \cite{Enqvist1978}, \cite{Jan1984}, \cite{Lee1982}, and~\cite{Nasari2012}. Recent theoretical developments include, e.g., \cite{Kato2019}, which conducts an in-depth investigation of the asymptotic behaviour of incomplete \mbox{U-statistics} in a high-dimensional framework, and~\cite{Lowe2021}, where a central limit theorem is obtained for triangular arrays. 
 Specific applications are treated in, e.g., \cite{Bertail06}, \cite{Bellet2016}, \cite{ClemProc2013}, \cite{Ment2016}, \cite{Papa2015}, and \cite{Rempala2003}. 
 
 The asymptotic distribution theory for incomplete \mbox{U-statistics} was first derived in~\cite{Jan1984}. However, like in subsequent papers, this was done there under the second-order moment assumption~${\rm E}[h^2(X_1,\ldots,X_\ell)]<\infty$. To the best of our knowledge, a consistency result that would only require the milder first-order moment assumption in~({\ref{firstmomentassumption}) is missing in the literature. In this note, we therefore establish such a consistency result (in the $L_1$ sense) for the incomplete \mbox{U-statistic} in~(\ref{UnN}). Our results are actually more general: for any~$p\geq 1$, we will obtain an 
 upper bound on the $L_1$ distance~${\rm E}[|U_{n,N}-\theta|]$ under the $p$th-order moment assumption stating that~${\rm E}[|h(X_1,\ldots,X_\ell)|^p]<\infty$. This will allow us to show consistency in the~$L_1$ sense and to derive an upper bound on the corresponding consistency rate. While the asymptotic distribution theory of incomplete $U$-statistics was derived under the assumption that~$n/N\to 0$ as~$n\to\infty$, we will make no restriction on the rate at which~$N$ diverges to infinity as~$n$ does. Actually, our upper bound on the $L_1$ distance~${\rm E}[|U_{n,N}-\theta|]$ even allows us to keep~$N$ fixed as~$n$ diverges to infinity.  

The outline of the paper is as follows. In Section~\ref{secBound}, we derive the aforementioned 
bound on~${\rm E}[|U_{n,N}-\theta|]$ under the $p$th-order moment assumptions. In Section~\ref{secRates}, we exploit the bound to establish consistency in the $L_1$ sense and to obtain, for each fixed~$p$, an upper bound on the corresponding consistency rate. In Section~\ref{secSampling}, we explain why our results actually extend to other classes of incomplete \mbox{U-statistics}. Finally, we perform a Monte-Carlo exercise in Section~\ref{secSimu} to illustrate our consistency results.


\section{A bound on the $L_1$ distance between~$U_{n,N}$ and~$\theta$}
\label{secBound}
\vspace{2mm}

Before stating the main result of this section, we provide the following comments on the random variables~$\alpha(i_1,\ldots,i_\ell)$ in~(\ref{UnN}), which will play an important role in our derivations. Under the considered sampling scheme, namely sampling without replacement, every fixed term from the complete \mbox{U-statistic}~$U_n$ in~(\ref{Un}) is selected with probability
$$
p_n
:=
{\rm E}[\alpha_n(i_1,\ldots,i_{\ell})]
=
\frac{N}{{n\choose \ell}}
. 
$$
Moreover, one has
\begin{equation}
\label{Varalpha}
\big|
{\rm E}[
(\alpha_n(i_1,\ldots,i_{\ell})-p_n)
(\alpha_n(j_1,\ldots,j_{\ell})-p_n)
]
\big|
\leq 
\Bigg\{
\begin{array}{ll}
\\[-5mm]
p_n & \textrm{if }(i_1,\ldots,i_{\ell})=(j_1,\ldots,j_{\ell}) 
\\[2mm]
p_n/{n\choose \ell} & \textrm{otherwise;}	
\end{array}
\end{equation}
see the proof of Lemma~2.1 in~\cite{Jan1984}. It is of course possible to provide exact values in~(\ref{Varalpha}), but our derivations below will only make use of the bounds in~(\ref{Varalpha}), which will allow us in Section~\ref{secSampling} to extend in an effortless way our results to other ways to sample terms from~$U_n$.

Under finite moment assumptions of order~$p$ (with~$p\in[1,\infty)$) on the kernel, we have the following bound on the $L_1$ distance between the incomplete statistic~$U_{n,N}$ in~(\ref{UnN}) and its expectation~$\theta$.  
\vspace{2mm}

\begin{theorem}
	\label{TheorMain}
Let~$X_1,X_2,\ldots$ be mutually independent and have a common distribution for which ${\rm E}[|h(X_1,\ldots,X_\ell)|^p]<\infty$ for some~$p\geq 1$.  Let~$(c_n)$ be a sequence that diverges to infinity. Then, 
\begin{equation}
\label{prerate}
{\rm E}[|U_{n,N}-\theta|]
=
O\Bigg(\frac{c_n^{(2-p)_+/2}}{\sqrt{\min(n,N)}}\Bigg)
+
o\bigg(\frac{1}{c_n^{p-1}}\bigg)
\end{equation}
as~$n\to \infty$, with~$\theta={\rm E}[h(X_1,\ldots,X_\ell)]$ 
$($throughout, we write~$x_+:=\max(x,0))$. 
\end{theorem}
\vspace{2mm}

The proof requires the following preliminary result, which we prove for the sake of completeness. 
\vspace{2mm}

\begin{lemma}
\label{LemMoments}
Let~$Z$ be a random variable taking values in~$\mathbb{R}^+$ and satisfying~${\rm E}[Z^p]<\infty$ for some~$p\geq 1$. Then,
$$
(i)
\ \
{\rm E}[ Z \mathbb{I}[Z>c] ]
=
o\Big(\frac{1}{c^{p-1}}\Big)
\quad
\textrm{and}
\quad
(i)
\ \
{\rm E}[ Z^2 \mathbb{I}[Z\leq c] ]
=
O(c^{(2-p)_+})
$$
as~$c\to \infty$.
\end{lemma}
\vspace{2mm}

\noindent
{\sc Proof of Lemma~\ref{LemMoments}.}
(i) Since the result trivially follows from Lebesgue's Dominated Convergence Theorem for~$p=1$, we may assume that~$p>1$. Letting~$q=p/(p-1)$ be the conjugate exponent of~$p$, H\"{o}lder's inequality yields
$
{\rm E}[ Z \mathbb{I}[Z>c] ]
\leq
({\rm E}[ Z^p ])^{1/p}
(P[X>c])^{1/q}
=
({\rm E}[ Z^p ])^{1/p}
(c^p P[X>c])^{1/q} c^{-p/q}
.
$
Since~$p/q=p-1$, the result follows from the fact that the moment assumption on~$Z$ ensures that~$c^p P[X>c]=o(1)$ as~$c\to \infty$ (see, e.g., the corollary on page~47 of~\cite{Ser1980}).
(ii) Let us assume that~$p\in[1,2)$ (the result is trivial for~$p\geq 2$). Then,
$
{\rm E}[ Z^2 \mathbb{I}[Z\leq c] ]
\leq
{\rm E}[ c^{2-p} Z^p \mathbb{I}[Z\leq c] ]
\leq
c^{2-p} {\rm E}[ Z^p ]
=
O(c^{2-p})
$  as~$c\to \infty$,
which establishes the result. 
\cqfd
\vspace{1mm}

We can now prove Theorem~\ref{TheorMain}. 
\vspace{3mm}

\noindent
{\sc Proof of Theorem~\ref{TheorMain}}
In this proof, $\mathcal{I}_{n,\ell}=\{I=(i_1,\ldots,i_\ell):1\leq i_1<\ldots<i_\ell\leq n\}$ will denote the collection of possible multi-indices~$I=(i_1,\ldots,i_\ell)$ in the complete \mbox{U-statistic}~$U_n$. This allows us to write
$$
U_{n}
=
\frac{1}{{n \choose \ell}}
\sum_{I\in \mathcal{I}_{n,\ell}}
h_n(X_I)
\ \ \
\textrm{and}
\ \ \
U_{n,N}
=
\frac{1}{N}
\sum_{I\in \mathcal{I}_{n,\ell}}
\alpha_n(I)
h_n(X_I)
,
$$
where~$X_I$ obviously stands for~$(X_{i_1},\ldots,X_{i_\ell})$. 
Based on the truncated kernel defined by
$$
h_n(x_1,\ldots,x_\ell)
:=
h(x_1,\ldots,x_\ell) 
\mathbb{I}[ |h(x_1,\ldots,x_\ell)| \leq c_n ]
,
$$
define then~$
\theta_n
:=
{\rm E}[
h_n(X_1,\ldots,X_\ell)
]
$,
$$
\tilde{U}_{n,N}
:=
\frac{1}{N}
\sum_{I\in \mathcal{I}_{n,\ell}}
\alpha_n(I)
h_n(X_I)
\ \ \
\textrm{and}
\ \ \
\bar{U}_{n,N}
:=
\frac{1}{N}
\sum_{I\in \mathcal{I}_{n,\ell}}
\alpha_n(I)
{\rm E}[h_n(X_I)]
=
\frac{\theta_n}{N}
\sum_{I\in \mathcal{I}_{n,\ell}}
\alpha_n(I)
.
$$
With this notation, write
\begin{eqnarray*}
{\rm E}[|U_{n,N}-\theta|]
\!&\!\! \leq \!\!&\! 
{\rm E}[|U_{n,N}-\tilde{U}_{n,N}|]
+
{\rm E}[|\tilde{U}_{n,N}-\bar{U}_{n,N}|]
+
{\rm E}[|\bar{U}_{n,N}-\theta_n|]
+
|\theta_n-\theta|
\\[2mm]
\!&\!\! =: \!\!&\! 
T_{1n}
+
T_{2n}
+ 
T_{3n}
+
T_{4n}
,
\end{eqnarray*}
say. We establish the result by proving that, for any~$j=1,2,3,4$, 
\begin{equation}
\label{ToshowTj}
T_{jn}
=
O\Bigg(\frac{c_n^{(2-p)_+/2}}{\sqrt{\min(n,N)}}\Bigg)
+
o\bigg(\frac{1}{c_n^{p-1}}\bigg)
\end{equation}
(throughout the proof, all~$o$'s and~$O$'s are as~$n\to\infty$).  
In view of the moment assumption on~$h(X_1,\ldots,X_\ell)$, Lemma~\ref{LemMoments}(i) yields that
\begin{eqnarray*}
T_{4n}
=
|\theta_n-\theta|
\!&\!\! \leq \!\!&\! 
{\rm E}[
|h_n(X_1,\ldots,X_\ell)-h(X_1,\ldots,X_\ell)|
]
\nonumber
\\[3mm]
\!&\!\! = \!\!&\! 
{\rm E}[
|h(X_1,\ldots,X_\ell)|
\mathbb{I}[|h(X_1,\ldots,X_\ell)|>c_n]
]
\nonumber
\\[2mm]
\!&\!\! = \!\!&\! 
o\bigg(\frac{1}{c_n^{p-1}}\bigg)
,
\label{Term1}
\end{eqnarray*}
so that~$T_{4n}$ satisfies~(\ref{ToshowTj}). 
By construction, 
\begin{equation}
\label{partT3n}
\bar{U}_{n,N}
=
\frac{\theta_n}{N}
\sum_{I\in \mathcal{I}_{n,\ell}}
\alpha_n(I)
=
\theta_n
,
\end{equation}
so that~$T_{3n}$ trivially satisfies~(\ref{ToshowTj}), too (the motivation to consider~$T_{3n}$ will be made clear when extending the results to other sampling schemes in Section~\ref{secSampling}). Now,
\begin{eqnarray*}
T_{1n}
\!&\!\!=\!\!&\!
{\rm E}[|U_{n,N}-\tilde{U}_{n,N}|]
\\[3mm]
\!&\!\!\leq \!\!&\!
\frac{1}{N}
\sum_{I\in \mathcal{I}_{n,\ell}}
{\rm E}[\alpha_n(I)]
{\rm E}[|h_n(X_I)-h(X_I)|]
\nonumber
\\[2.5mm]
\!&\!\!=\!\!&\!
{\rm E}[
|h(X_1,\ldots,X_\ell)|
\mathbb{I}[|h(X_1,\ldots,X_\ell)|> c_n]
]
\nonumber
\\[2.5mm]
\!&\!\!=\!\!&\!
o\bigg(\frac{1}{c_n^{p-1}}\bigg)
,
\label{Term2}
\end{eqnarray*}
so that~$T_{1n}$ satisfies~(\ref{ToshowTj}) as well. It remains to treat~$T_{2n}$, which is more delicate. First note that
$$
\tilde{U}_{n,N} - \bar{U}_{n,N}
=
\frac{1}{N}
\sum_{I\in\mathcal{I}_{n,\ell}}
\alpha_n(I)
g_n(X_I)
,
$$
where~$g_n$ is the centered truncated kernel defined by
$$
g_n(x_1,\ldots,x_\ell)
:=
h_n(x_1,\ldots,x_\ell) 
-
{\rm E}[h_n(X_1,\ldots,X_\ell)]
.
$$
Since~$g_n$ is bounded by construction, we can consider
\begin{eqnarray*}
T_{2n}^2
\leq
{\rm E}[(\tilde{U}_{n,N} - \bar{U}_{n,N})^2]
\!&\!\! = \!\!&\! 
\frac{1}{N^2}
\sum_{I,J\in\mathcal{I}_{n,\ell}}
{\rm E}[
\alpha_n(I)
\alpha_n(J)
g_n(X_I)
g_n(X_J)
]
\\[2mm]
\!&\!\! = \!\!&\! 
\frac{1}{N^2}
\sum_{I,J\in\mathcal{I}_{n,\ell}}
{\rm E}[
\alpha_n(I)
\alpha_n(J)
]
{\rm E}[
g_n(X_I)
g_n(X_J)
]
\\[2mm]
\!&\!\! \leq \!\!&\! 
a_n
+
b_n
,
\end{eqnarray*}
with 
(recall~(\ref{Varalpha}))
$$
a_n
:=
\frac{p_n}{N^2}
\sum_{I\in\mathcal{I}_{n,\ell}}
{\rm E}[
g_n^2(X_I)
]
\quad
\textrm{and}
\quad
b_n
:=
\frac{2p_n^2}{N^2}
\sum_{I,J\in\mathcal{I}_{n,\ell}, I\neq J}
{\rm E}[
g_n(X_I)
g_n(X_J)
]
;
$$
for~$b_n$, we used~(\ref{Varalpha}) to obtain 
\begin{equation}
\label{heree}
{\rm E}[
\alpha_n(I)
\alpha_n(J)
]
=
{\rm Cov}[
\alpha_n(I)
,
\alpha_n(J)
]
+
p_n^2
\leq
\frac{p_n}{{n\choose \ell}}
+
p_n^2
=
\frac{p_n^2}{N}
+
p_n^2
\leq
2p_n^2
.
\end{equation}
For~$a_n$, Lemma~\ref{LemMoments}(ii) yields
\begin{eqnarray}
a_n
\!&\!\! = \!\!&\!
\frac{p_n}{N^2}
{n \choose \ell}
{\rm E}[
g_n^2(X_1,\ldots,X_\ell)
]
=
\frac{1}{N}
{\rm Var}[
h_n(X_1,\ldots,X_\ell)
]
\nonumber
\\[2.5mm]
\!&\!\!\leq\!\!&\!
\frac{1}{N}
{\rm E}[
h_n^2(X_1,\ldots,X_\ell)
]
=
\frac{1}{N}
{\rm E}[
h^2(X_1,\ldots,X_\ell)
\mathbb{I}[ |h(X_1,\ldots,X_\ell)| \leq c_n ]
]
\nonumber
\\[2.5mm]
\!&\!\! = \!\!&\!
O\bigg(\frac{c_n^{(2-p)_+}}{\min(n,N)}\bigg)
.
\label{T1}
\end{eqnarray}
Consider then~$b_n$. If the multi-indices~$I$ and~$J$ do not have at least one index in common, then mutual independence of the~$X_i$'s readily yields
$$
{\rm E}[
g_n(X_I)
g_n(X_J)
]
=
{\rm E}[g_n(X_I)]
{\rm E}[g_n(X_J)]
=
0
.
$$
Out of the 
$$
{n \choose \ell} 
\bigg( {n \choose \ell} - 1 \bigg)
$$
terms in~$b_n$, we may thus restrict to the 
$$
m_n
:=
{n \choose \ell} 
\bigg( {n \choose \ell} - 1 \bigg)
-
{n \choose \ell} 
{n-\ell \choose \ell} 
$$
terms for which~$I$ and~$J$ have at least one index in common. The Cauchy--Schwarz inequality and Lemma~\ref{LemMoments}(ii) then provide
\begin{eqnarray}
b_n
\!&\!\!\leq\!\!&\!
\frac{2m_n p_n^2}{N^2}
{\rm E}[
g_n^2(X_1,\ldots,X_\ell)
]
=
\frac{2m_n}{{n \choose \ell}^2}
{\rm Var}[h_n(X_1,\ldots,X_\ell)]
\nonumber
\\[2mm]
\!&\!\!\leq\!\!&\!
\frac{2m_n}{{n \choose \ell}^2}
{\rm E}[h_n^2(X_1,\ldots,X_\ell)]
\leq
\frac{2m_n}{{n \choose \ell}^2}
{\rm E}[
h^2(X_1,\ldots,X_\ell
\mathbb{I}[ |h(X_1,\ldots,X_\ell)| \leq c_n ]
]
\nonumber
\\[2mm]
\!&\!\!=\!\!&\!
\frac{2m_n}{{n \choose \ell}^2}
O(c_n^{(2-p)_+})
=
O\bigg(\frac{c_n^{(2-p)_+}}{\min(n,N)}\bigg)
,
\label{T2}
\end{eqnarray}
where we used the fact that
\begin{eqnarray*}
\lefteqn{
\frac{m_n}{{n \choose \ell}^2}
=
1
-
\frac{1}{{n \choose \ell}} 
-
\frac{{n-\ell \choose \ell}}{{n \choose \ell}}
}
\\[2mm]
& & 
\hspace{3mm}
=
1
-
\frac{(n-\ell)(n-\ell-1)\ldots(n-2\ell+1)}{n(n-1)\ldots(n-\ell+1)}
+
O\Big(\frac{1}{n}\Big)
\\[2mm]
& & 
\hspace{3mm}
=
1
-
({\textstyle{1-\frac{\ell}{n}}})
({\textstyle{1-\frac{\ell}{n-1}}})
\ldots
({\textstyle{1-\frac{\ell}{n-\ell+1}}})
+
O\Big(\frac{1}{n}\Big)
\\[2mm]
& & 
\hspace{3mm}
=
O\Big(\frac{1}{n}\Big)
.
\end{eqnarray*}
Thus, (\ref{T1})--(\ref{T2}) yield
$$
T_{2n}
=
O\Bigg(\frac{c_n^{(2-p)_+/2}}{\sqrt{\min(n,N)}}\Bigg)
.
$$
This shows that~$T_{2n}$ satisfies~(\ref{ToshowTj}), too, which establishes the result. 
\cqfd
\vspace{1mm}

Under the minimal assumption making~$\theta$ well-defined, namely under the assumption that~${\rm E}[|h(X_1,\ldots,X_\ell)|]<\infty$ ($p=1$), a direct corollary of Theorem~\ref{TheorMain} is that 
$$
{\rm E}[|U_{n,N}-\theta|]=o(1)
$$   
as soon as~$\min(n,N)\to\infty$ (this follows by taking, e.g., $c_n=\log(\min(n,N))$). This establishes convergence in the~$L_1$ sense (hence also in probability) of~$U_{n,N}$ to its expectation~$\theta$. Interestingly, convergence holds even when~$N$ would diverge to~$n$ arbitrarily slowly as a function of~$n$, whereas classical asymptotic theory of incomplete \mbox{U-statistics} requires that~$n/N\to 0$ as~$n\to\infty$; see~\cite{Jan1984}.


\section{Consistency rates}
\label{secRates}
\vspace{2mm}

As just pointed out, first-order moment assumptions on the kernel~$h(X_1,\ldots,X_\ell)$  are sufficient to ensure convergence of incomplete \mbox{U-statistics} to their expectation. Actually, the bound obtained in Theorem~\ref{TheorMain} also allows us to comment on the convergence rate that is achieved when assuming that~${\rm E}[|h(X_1,\ldots,X_\ell)|^p]<\infty$.

Consider first the case~$p\geq 2$. Taking~$c_n=\sqrt{\min(n,N)}$ in~(\ref{prerate}) then provides  
\begin{equation}
\label{rate2ndorder}
{\rm E}[|U_{n,N}-\theta|]
=
O\bigg(\frac{1}{\sqrt{\min(n,N)}}\bigg)
,
\end{equation}
which is compatible with the (optimal) root-$n$ consistency rate that was obtained in~\cite{Jan1984} under second-order moment assumptions (\cite{Jan1984} restricts to the case where~$n/N=o(1)$, in which case the optimal rate is indeed $1/\sqrt{n}$, but it should be clear that the optimal rate generalizes into the one in~(\ref{rate2ndorder}) when no assumption is made on the rate at which~$N$ diverges to infinity as~$n$ does). Now, consider the case~$p\in[1,2)$, that is the one of main interest in the present work. Would both terms in the righthand side of~(\ref{prerate}) be big~$O$'s, an upper bound on the consistency rate could be obtained by choosing~$c_n$ so that
$$
\frac{c_n^{(2-p)_+/2}}{\sqrt{\min(n,N)}}
\sim
\frac{1}{c_n^{p-1}}
$$
(where~$a_n\sim b_n$ means that $a_n=O(b_n)$ and $b_n=O(a_n)$), that is, by taking
$
c_n
\sim
\{\min(n,N)\}^{1/p}
$,
which would provide
$$
{\rm E}[|U_{n,N}-\theta|]
=
O\bigg(\frac{1}{\{\min(n,N)\}^{(p-1)/p}}\bigg)
.
$$
Since the second term in the righthand side of~(\ref{prerate}) is a little~$o$ rather than a big~$O$, it is possible to take~$c_n=r_n\{\min(n,N)\}^{1/p}$, with a suitable sequence~$(r_n)$ that is~$o(1)$, to obtain
\begin{equation}
\label{rate2ndorder}
{\rm E}[|U_{n,N}-\theta|]
=
o\bigg(\frac{1}{\{\min(n,N)\}^{(p-1)/p}}\bigg)
\end{equation}
for~$p\in[1,2)$. This is only an upper bound on the convergence rate in the $L_1$ sense, but the natural values of this upper bound for~$p=1$ and~$p\stackrel{<}{\to} 2$, namely
$$
{\rm E}[|U_{n,N}-\theta|]
=
o(1)
\ \ \ 
\textrm{ and }
\ \ \ 
{\rm E}[|U_{n,N}-\theta|]
=
o\bigg(\frac{1}{\sqrt{\min(n,N)}}\bigg)
,
$$
respectively, lead us to conjecture that this is actually the exact convergence rate. Note also that, as expected, the larger~$p\in[1,2)$, the better the convergence rate in~(\ref{rate2ndorder}).



\section{Extensions to other sampling schemes}
\label{secSampling}
\vspace{2mm}

The incomplete \mbox{U-statistic}~$U_{n,N}$ considered in the earlier sections is the one in~(\ref{UnN}), that is obtained by \emph{sampling without replacement} $N$ terms from the~${n\choose \ell}$ terms in the complete \mbox{U-statistic}~$U_n$ in~(\ref{Un}). Other classical sampling schemes are possible, including both following ones:
\vspace{1mm}

\begin{itemize}
\item \emph{Sampling with replacement:} in this scheme, the~${n\choose \ell}$ random variables $\alpha(i_1,\ldots,i_\ell)$, with values in~$\{0,1,\ldots,N\}$, are the marginals of a multinormal random vector with count parameter~$N$ and  homogeneous probabilities~$1/{n\choose \ell},\ldots,1/{n\choose \ell}$ on each of the~${n\choose \ell}$ possible outcomes. Here, $N$ is an arbitrary positive integer and we let~$p_n=N/{n\choose \ell}$.
\vspace{2mm}
\item \emph{Bernoulli sampling:} in this scheme, the~${n\choose \ell}$ random variables $\alpha(i_1,\ldots,i_\ell)$ form a random sample from the Bernoulli distribution with success probability~$p_n=N/{n\choose \ell}$, with~$N\in(0,{n\choose \ell}]$. 
\end{itemize}
For \emph{sampling with replacement}, it is easy to check that, using the notation from the proof of Theorem~\ref{TheorMain},  we still have
\vspace{3mm}
\begin{equation}
\label{AA1}
N\geq 1
,
\vspace{3mm}
\end{equation}
\begin{equation}
\label{AA2}
\sum_{I\in\mathcal{I}_{n,\ell}}
\alpha_n(I)
=
N
,
\end{equation}
\begin{equation}
\label{AA3}
{\rm E}[\alpha_n(I)]
=
p_n
=
\frac{N}{{n\choose \ell}}
\end{equation}
and that the covariance bounds in~(\ref{Varalpha}) still hold. Since these were the only assumptions used on~$N$ and on the random variables~$\alpha(i_1,\ldots,i_\ell)$ in the proof of Theorem~\ref{TheorMain}, we conclude that this result, hence also its consequences described at the end of Section~\ref{secBound} and in Section~\ref{secRates}, extend to incomplete \mbox{U-statistics} that are obtained from \emph{sampling with replacement}.

The situation is slightly more complicated for \emph{Bernoulli sampling}, for which~(\ref{AA3}) and the covariance bounds in~(\ref{Varalpha}) still hold, but for which~(\ref{AA1})--(\ref{AA2}) do not hold in general. An inspection of the proof of Theorem~\ref{TheorMain} reveals that the proof can easily be adapted to the case where~$N$ remains away from zero --- which is obviously an extremely mild assumption since consistency requires anyway that~$N$ diverges to infinity as~$n$ does. The fact that~(\ref{AA2}) that does not hold for \emph{Bernoulli sampling} is more serious, since it implies that~(\ref{partT3n}) does not hold for this sampling scheme, which forces us to control~$T_{3n}$: since
$$
\bar{U}_{n,N}-\theta_n
=
\frac{\theta_n}{N}
\sum_{I\in \mathcal{I}_{n,\ell}}
(\alpha_n(I)-p_n)
,
$$
using~(\ref{Varalpha}) and the fact that~$\theta_n=\theta+o(1)$ in the proof of Theorem~\ref{TheorMain} yields
\begin{eqnarray*}
T_{3n}^2
\leq
{\rm E}[ (\bar{U}_{n,N}-\theta_n)^2 ]
\!&\!\! = \!\!&\! 
\frac{\theta_n^2}{N^2}
\sum_{I,J\in \mathcal{I}_{n,\ell}}
\!\!
{\rm E}[
(\alpha_n(I)-p_n)
(\alpha_n(J)-p_n)
]
\nonumber
\\[2mm]
\!&\!\! \leq \!\!&\! 
\frac{\theta_n^2}{N^2}
\bigg\{
{n \choose \ell}
\times
p_n
+
{n \choose \ell} \bigg({n \choose \ell}-1\bigg)
\times
\frac{p_n}{{n \choose \ell}}
\bigg\}
\nonumber
\\[2mm]
\!&\!\! \leq \!\!&\! 
\frac{2\theta_n^2 {n \choose \ell}p_n}{N^2}
=
\frac{2\theta_n^2}{N}
=
O\Big(\frac{1}{N}\Big)
\end{eqnarray*}
as~$n$ diverges to infinity.
Thus, 
$$
T_{3n}
=
O\bigg(\frac{1}{\sqrt{N}}\bigg)
=
O\Bigg(\frac{c_n^{(2-p)_+/2}}{\sqrt{\min(n,N)}}\Bigg)
$$
as~$n$ diverges to infinity, which shows that~$T_{3n}$ satisfies~(\ref{ToshowTj}) for \emph{Bernoulli sampling}, too. Therefore, we conclude that all consistency results of the paper apply to this third sampling scheme as well.


\section{A numerical illustration}
\label{secSimu}
\vspace{2mm}

We end this paper with a Monte-Carlo exercise we performed in order to illustrate our consistency results. In this numerical exercise, we focus on scalar observations ($\mathcal{X}=\mathbb{R}$) and on the kernel~$h:\mathbb{R}^2\to \mathbb{R}$ defined by
$$
h(x_1,x_2) = |x_1-x_2|
.
$$
For each~$n\in\{50,100,150,\ldots,400\}$, we generated $M=12,000$ random samples of size~$n$ from the Student~$t$ distribution with~$\nu=1.5$ (for which the kernel~$h(X_1,X_2)$ has infinite second-order moments). In each of the resulting samples, we evaluated the complete \mbox{U-statistic}
\begin{equation}
\label{completesimu}
U_n 
=
\frac{1}{{n\choose 2}} 
\sum_{1\leq i<j\leq n} 
|X_i-X_j|
\end{equation}
and nine incomplete \mbox{U-statistics} of the form
\begin{equation}
\label{incompletesimu}
U_{n,N} 
=
\frac{1}{N} 
\sum_{1\leq i<j\leq n} 
\alpha(i,j)
|X_i-X_j|
,
\end{equation}
obtained by combining one of the three possible sampling schemes (\emph{sampling without replacement}, \emph{sampling with replacement}, and \emph{Bernoulli sampling}) and a value of~$N\in\{n^{3/2},n,n^{2/3}\}$ (recall that the complete statistic involves~${n\choose 2}=n(n-1)/2$ terms). Complete and incomplete \mbox{U-statistics} here are consistent for~$\theta={\rm E}[|X_1-X_2|]
$. For each value of~$n$ and for each given complete or incomplete \mbox{U-statistic},~$V_n$ say, we evaluated the empirical $L_1$ distance
\begin{equation}
\label{empirL1dist}
\frac{1}{M}
\sum_{m=1}^M
|V_n(m)-\theta|
\end{equation}
between~$V_n$ and~$\theta$, where~$V_n(m)$ denotes the value taken by~$V_n$ in the $m$th random sample of size~$n$ that was generated above. In Figure~\ref{Fig1}, we plot the quantities in~(\ref{empirL1dist}) versus the sample size~$n$ for the complete \mbox{U-statistic} and each of the nine incomplete statistics mentioned above. To explore the sensitivity of the results to moment assumptions on~$h(X_1,X_2)$, we repeated this exercise for random samples generated from Student~$t$ distributions with~$\nu=1.8$, $2.1$ and~$4.1$ degrees of freedom (note that the kernel~$h(X_1,X_2)$ has infinite second-order moments if and only if~$\nu\leq 2$).

\begin{figure}[h!]
	\includegraphics[width=\textwidth]{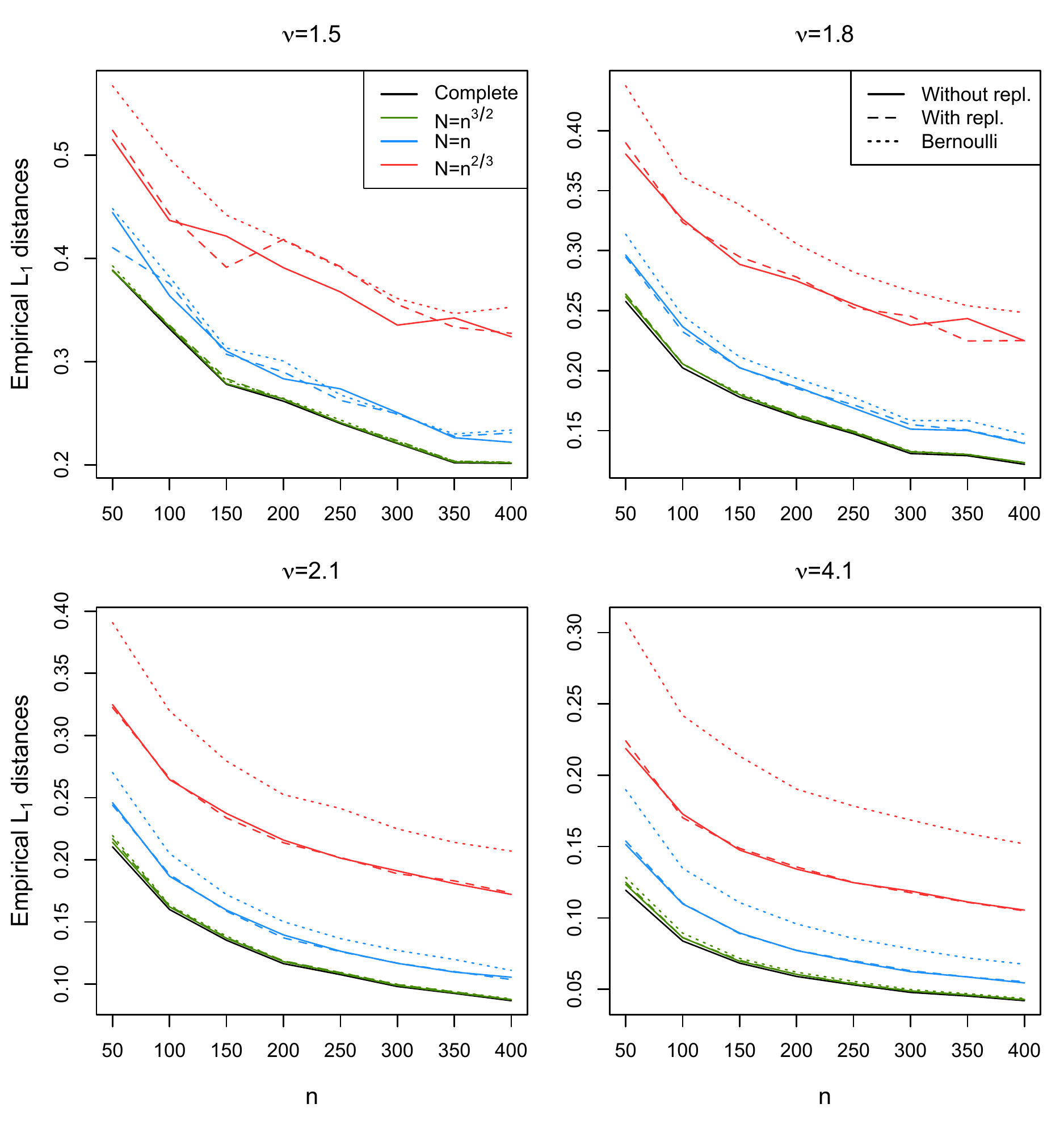}
\caption{Empirical $L_1$ distances in~(\ref{empirL1dist}) versus the sample size~$n$ for the complete \mbox{U-statistic} in~(\ref{completesimu}) and the nine versions of the incomplete \mbox{U-statistic} in~(\ref{completesimu}), when observations are randomly sampled from the Student~$t$ distribution 
with~$\nu=1.5$ degrees of freedom (top left),
$\nu=1.8$ degrees of freedom (top right),
$\nu=2.1$ degrees of freedom (bottom left),
or 
$\nu=4.1$ degrees of freedom
(bottom right).
The nine version of the incomplete \mbox{U-statistic} result from the combination of one of the three possible sampling schemes (\emph{sampling without replacement}, \emph{sampling with replacement}, and \emph{Bernoulli sampling}) and one of the three considered values of~$N\in\{n^{3/2},n,n^{2/3}\}$.}
\label{Fig1}
\end{figure}

Figure~\ref{Fig1} clearly supports our theoretical results. Irrespective of moment assumptions, empirical $L_1$ distances of all complete and incomplete statistics decrease to zero as~$n$ increases. As expected, the smaller~$N$, the larger these distances, but the results also reveal that, irrespective of the adopted sampling scheme, $N=n^{3/2}$ provide virtually the same $L_1$ distances as the complete U-statistics. Maybe less importantly, this numerical exercise also shows that, while sampling with replacement and sampling without replacement provide the same performances, Bernoulli sampling yields higher $L_1$ distances, with a difference that increases as~$N$ decreases.


\bibliographystyle{plain}
\bibliography{paper}

\end{document}